\documentclass{amsart}

\usepackage{graphicx}
\usepackage{tikz}

\usepackage{amsmath,amsfonts,amsthm,amssymb,amscd}
\usepackage{hyperref}

\usepackage[alphabetic]{amsrefs}

\newcommand*{\textlabel}[2]{%
  \edef\@currentlabel{#1}
  \phantomsection
  #1\label{#2}
}
\makeatother

\newtheorem*{thmspec}{Theorem}

\theoremstyle{definition}

\numberwithin{equation}{section}

\begin{document}

\title[Hex implies Y]{Hex implies Y}

\author[T. Prytu{\l}a]{Tomasz Prytu{\l}a}

 \thanks{ The author was supported by the EU Horizon 2020 program under the Marie Sk{\l}odowska-Curie grant agreement no. 713683 (COFUNDfellowsDTU)}

\address{Department of applied mathematics and computer science, Technical University of Denmark, Lyngby, Denmark}

\email{tompr@dtu.dk}



\begin{abstract}
	We give a simple and short proof of the fact that the board game of Y cannot end in a draw. Our proof, based on the analogous result for the game of Hex (the so-called ‘Hex Theorem’), is purely topological and does not depend on the shape of the board. We also include a simplified version of Gale's proof of Hex Theorem.
\end{abstract}

\maketitle

\section{Introduction}

\begin{figure}[!h]
  \centering
    \begin{tikzpicture}[scale=0.25]

      \definecolor{vlgray}{RGB}{230,230,230}
      \definecolor{red}{RGB}{255,000,000}

      \definecolor{blue}{RGB}{000,000,255}

\begin{scope}[shift={(-12,0)}]

\draw[fill=red, red] (3,0.5)   circle [radius=0.65];

\draw[fill=red, red] (2,2.25)   circle [radius=0.65];

\draw[fill=red, red] (3,4)   circle [radius=0.65];

\draw[fill=red, red] (4,5.75)   circle [radius=0.65];

\draw[fill=red, red] (3,7.5)   circle [radius=0.65];

\draw[fill=red, red] (4,9.25)   circle [radius=0.65];

\draw[fill=red, blue] (-1,4)   circle [radius=0.65];

\draw[fill=blue, blue] (1,4)   circle [radius=0.65];

\draw[fill=blue, blue] (5,4)   circle [radius=0.65];

\draw[fill=blue, blue] (6,5.75)   circle [radius=0.65];

\draw[fill=blue, blue] (8,5.75)   circle [radius=0.65];

\foreach \j in {0,...,5}{
	\foreach \i in {0,...,5}{
    
    	\begin{scope}[shift={(2*\i-\j,\j*1.75)}]
      		\draw (0,0)--(0,1)--(1,1.75)--(2,1)--(2,0)--(1,-0.75)--(0,0);
      	\end{scope}
    }
}

\begin{scope}

\draw[ blue, very thick] (1,-0.75)--(0,0)--(0,1)--(-1, 1.75)--(-1,2.75)--(-2, 3.5)--(-2,4.5)--(-3,5.25)--(-3,6.25)--(-4,7)--(-4,8)--(-5,8.75)--(-5,9.75);

\end{scope}

\begin{scope}[shift={(12,0)}]

\draw[ blue, very thick] (0,0)--(0,1)--(-1, 1.75)--(-1,2.75)--(-2, 3.5)--(-2,4.5)--(-3,5.25)--(-3,6.25)--(-4,7)--(-4,8)--(-5,8.75)--(-5,9.75)--(-6, 10.5);
\end{scope}

\begin{scope}[shift={(0,0)}]

\draw[ red, very thick] (1,-0.75)--(2,0)--(3,-0.75)--(4, 0)--(5,-0.75)--(6, 0)--(7,-0.75)--(8,0)--(9,-0.75)--(10,0)--(11,-0.75)--(12,0);
\end{scope}

\begin{scope}[shift={(-6,10.5)}]

\draw[ red, very thick] (1,-0.75)--(2,0)--(3,-0.75)--(4, 0)--(5,-0.75)--(6, 0)--(7,-0.75)--(8,0)--(9,-0.75)--(10,0)--(11,-0.75)--(12,0);
\end{scope}

\end{scope}


  \begin{scope}[shift={(12,8.75)}]

 \begin{scope}[shift={(-3,-8.75)}]

\draw[fill=red, red] (3,0.5)   circle [radius=0.65];

\draw[fill=red, red] (2,2.25)   circle [radius=0.65];

\draw[fill=red, red] (3,4)   circle [radius=0.65];

\draw[fill=red, red] (4,5.75)   circle [radius=0.65];

\draw[fill=red, red] (3,7.5)   circle [radius=0.65];

\draw[fill=red, red] (6,5.75)   circle [radius=0.65];

\draw[fill=blue, blue] (1,4)   circle [radius=0.65];

\draw[fill=blue, blue] (4,2.25)   circle [radius=0.65];

\draw[fill=blue, blue] (6,2.25)   circle [radius=0.65];

\draw[fill=blue, blue] (7,4)   circle [radius=0.65];

\draw[fill=blue, blue] (7,0.5)   circle [radius=0.65];
\end{scope}
     
\foreach \j in {0,...,5}{
	\foreach \i in {0,...,\j}{
    
    	\begin{scope}[shift={(2*\i-\j,-\j*1.75)}]
      		\draw (0,0)--(0,1)--(1,1.75)--(2,1)--(2,0)--(1,-0.75)--(0,0);
      	\end{scope}
    }
}

\begin{scope}[shift={(-5,-8.75)},xscale=-1]

\draw[thick] (0,0)--(0,1)--(-1, 1.75)--(-1,2.75)--(-2, 3.5)--(-2,4.5)--(-3,5.25)--(-3,6.25)--(-4,7)--(-4,8)--(-5,8.75)--(-5,9.75)--(-6, 10.5);
\end{scope}

\begin{scope}[shift={(7,-8.75)}]

\draw[ thick] (0,0)--(0,1)--(-1, 1.75)--(-1,2.75)--(-2, 3.5)--(-2,4.5)--(-3,5.25)--(-3,6.25)--(-4,7)--(-4,8)--(-5,8.75)--(-5,9.75)--(-6, 10.5);

\end{scope}
\begin{scope}[shift={(-5,-8.75)}]
\draw[thick] (0,0)--(1,-0.75)--(2,0)--(3,-0.75)--(4, 0)--(5,-0.75)--(6, 0)--(7,-0.75)--(8,0)--(9,-0.75)--(10,0)--(11,-0.75)--(12,0);


\end{scope}

\end{scope}

    \end{tikzpicture}
  \caption{Boards for games Hex and Y. On both boards the position is won by red.}
  \label{fig:gamescells}
\end{figure}
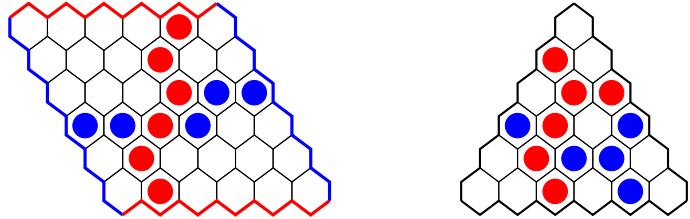

  Games of Hex and Y are turn-based, two player, abstract strategy games, belonging to a family of connection games.  
  Both games, in their original form, are played on a subset of a hexagonal tessellation of the plane, where the players (commonly denoted by red and blue) take turns to place stones on unoccupied hexagonal cells. Once the stone is placed it cannot be moved or removed. Hex is played on a rectangular board, whose two pairs of opposite sides are denoted by red, and blue respectively. The goal for the red player is to create a path of red stones joining two red sides, and the goal for the blue player is to create a blue path joining two blue sides.
  The game of Y is played on a triangular board, and the goal for both players is to create a connected chain of stones joining all three sides of the board (such a chain will generically have a shape of the letter Y, hence the name of the game). Figure~\ref{fig:gamescells} shows exemplary boards and positions for Hex and Y.

  Hex was first described by Hein in 1942, and rediscovered by Nash in 1948, and it has been quite extensively studied since (see \cite{Gardner59} for a brief history of Hex). The game of Y is an interesting variation on the theme, and was discovered independently by Milnor and Shannon in 1950s, and rediscovered by Schensted and Titus in 1953 \cites{Gardner08, Nash, SchenTit}. An appealing feature of Y is that the objective of the game for both players is the same, as opposed to Hex.

  One of the main mathematical properties of Hex is that it cannot end in a draw. More precisely, once the Hex board is completely filled with stones, then there is either a red path joining the red sides, or a blue path joining the blue sides, but not both. This is known as ‘Hex Theorem’. The same is true for the game of Y; we will refer to it as ‘Y Theorem’. Over the years multiple proofs of Hex Theorem appeared (see e.g.,\ \cites{Beck, Gale, Berman}), however, there are few arguments for Y Theorem, and all the proofs known to us are based on the combinatorial properties of the board and use some sort of induction or recursion \cites{Rijswick,HayRijs}. The main purpose of this note is to show that Y~Theorem is in fact equivalent to Hex Theorem, and thus obtain a new, very simple and purely topological proof of Y Theorem.

  \begin{thmspec} 
    Hex Theorem and Y Theorem are equivalent.
  \end{thmspec}

  The fact that Y Theorem implies Hex Theorem is an easy trick, as a game of Hex can be seen as a continuation of a game of Y from a certain position. This was observed by Schensted \cite{SchenTit}, \cite[Section 4.7]{Rijswick}. Our efforts thus go into proving the converse implication. The key idea is to ‘double’ the Y board by reflecting it along one of its sides, and then treat the resulting $4$--gon as a Hex board.

  We present both theorems in a generalized form, where the game is played on an arbitrary triangulation of a $2$--disk. This generality, besides the obvious benefit of obtaining a more general result, also allows us to use precise language of graph theory.

  For the sake of completeness we also include the proof of Hex Theorem. We claim no originality for the ideas behind this proof, apart from a minor simplification at the very end, which allows us to keep the proof concise. This note may thus also serve as a short and self-contained proof of Hex Theorem in a generalized setting.

\section{Hex and Y in a generalized form}

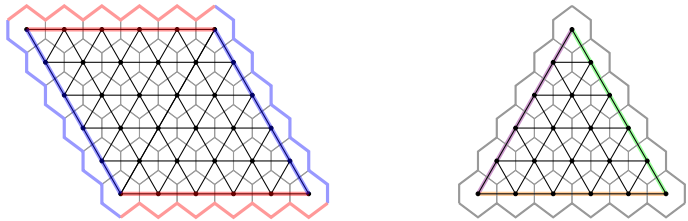
\begin{figure}[!h]
  \centering
    \begin{tikzpicture}[scale=0.25]

      \definecolor{vlgray}{RGB}{230,230,230}
      \definecolor{red}{RGB}{255,000,000}

      \definecolor{blue}{RGB}{000,000,255}

\begin{scope}[shift={(-12,0)}]


\foreach \j in {0,...,5}{
	\foreach \i in {0,...,5}{
    
    	\begin{scope}[shift={(2*\i-\j,\j*1.75)}]
      		\draw[black!40] (0,0)--(0,1)--(1,1.75)--(2,1)--(2,0)--(1,-0.75)--(0,0);

			\draw[fill=black, black] (1,0.5)   circle [radius=0.1];
      	\end{scope}
    }
}

\foreach \j in {0,...,5}{

    	\begin{scope}[shift={(-\j,\j*1.75)}]
      		\draw (1,0.5)--(11,0.5);

      	\end{scope}

      	  	\begin{scope}[shift={(2*\j,0)}]
      		\draw (1,0.5)--(-4,9.25);
      		
      	\end{scope}
}

\foreach \i in {0,...,5}{ 

\begin{scope}[shift={(10-2*\i,0)}]

      		\draw (1,0.5)--(1+ \i,0.5+ 1.75*\i);

\end{scope}

}

\foreach \i in {0,...,5}{ 

\begin{scope}[shift={(-5+2*\i,8.75)}]

      		\draw (1,0.5)--(1-\i,0.5- 1.75*\i);

\end{scope}

}

\begin{scope}

\draw[ blue!40, very thick] (1,-0.75)--(0,0)--(0,1)--(-1, 1.75)--(-1,2.75)--(-2, 3.5)--(-2,4.5)--(-3,5.25)--(-3,6.25)--(-4,7)--(-4,8)--(-5,8.75)--(-5,9.75);

\end{scope}

\begin{scope}[shift={(12,0)}]

\draw[ blue!40, very thick] (0,0)--(0,1)--(-1, 1.75)--(-1,2.75)--(-2, 3.5)--(-2,4.5)--(-3,5.25)--(-3,6.25)--(-4,7)--(-4,8)--(-5,8.75)--(-5,9.75)--(-6, 10.5);
\end{scope}

\begin{scope}[shift={(0,0)}]

\draw[ red!40, very thick] (1,-0.75)--(2,0)--(3,-0.75)--(4, 0)--(5,-0.75)--(6, 0)--(7,-0.75)--(8,0)--(9,-0.75)--(10,0)--(11,-0.75)--(12,0);
\end{scope}

\begin{scope}[shift={(-6,10.5)}]

\draw[ red!40, very thick] (1,-0.75)--(2,0)--(3,-0.75)--(4, 0)--(5,-0.75)--(6, 0)--(7,-0.75)--(8,0)--(9,-0.75)--(10,0)--(11,-0.75)--(12,0);
\end{scope}

\draw[ red, opacity=0.4, ultra thick] (1,0.5)--(11,0.5);

\draw[ red, opacity=0.4, ultra thick] (-4,9.25)--(6,9.25);

\draw[ blue, opacity=0.4, ultra thick] (1,0.5)--(-4,9.25);

\draw[ blue, opacity=0.4, ultra thick] (6,9.25)--(11,0.5);

\end{scope}


  \begin{scope}[shift={(12,8.75)}]

\draw[ green, opacity=0.4, ultra thick] (1,0.5)--(6,-8.25);

\draw[ orange, opacity=0.4, ultra thick] (-4,-8.25)--(6,-8.25);

\draw[ violet, opacity=0.4, ultra thick] (1,0.5)--(-4,-8.25);

\foreach \j in {0,...,5}{
	\foreach \i in {0,...,\j}{
    
    	\begin{scope}[shift={(2*\i-\j,-\j*1.75)}]
      		\draw[black!40] (0,0)--(0,1)--(1,1.75)--(2,1)--(2,0)--(1,-0.75)--(0,0);
      		\draw[fill=black, black] (1,0.5)   circle [radius=0.1];

      	\end{scope}
    }
}

\draw (1,0.5)--(6,-8.25);
\draw (0,-1.25)--(4,-8.25);
\draw (-1,-3)--(2,-8.25);
\draw (-2,-4.75)--(0,-8.25);
\draw (-3,-6.5)--(-2,-8.25);

\draw (1,0.5)--(-4,-8.25);   
\draw (2,-1.25)--(-2,-8.25);   
\draw (3,-3)--(0,-8.25);   
\draw (4,-4.75)--(2,-8.25);   
\draw (5,-6.5)--(4,-8.25);

\draw (-4,-8.25)--(6,-8.25);
\draw (-3,-6.5)--(5,-6.5);
\draw (-2,-4.75)--(4,-4.75);
\draw (-1,-3)--(3,-3);
\draw (0,-1.25)--(2,-1.25);

\begin{scope}[shift={(-5,-8.75)},xscale=-1]

\draw[black!40, thick] (0,0)--(0,1)--(-1, 1.75)--(-1,2.75)--(-2, 3.5)--(-2,4.5)--(-3,5.25)--(-3,6.25)--(-4,7)--(-4,8)--(-5,8.75)--(-5,9.75)--(-6, 10.5);
\end{scope}

\begin{scope}[shift={(7,-8.75)}]

\draw[black!40, thick] (0,0)--(0,1)--(-1, 1.75)--(-1,2.75)--(-2, 3.5)--(-2,4.5)--(-3,5.25)--(-3,6.25)--(-4,7)--(-4,8)--(-5,8.75)--(-5,9.75)--(-6, 10.5);

\end{scope}
\begin{scope}[shift={(-5,-8.75)}]
\draw[black!40, thick] (0,0)--(1,-0.75)--(2,0)--(3,-0.75)--(4, 0)--(5,-0.75)--(6, 0)--(7,-0.75)--(8,0)--(9,-0.75)--(10,0)--(11,-0.75)--(12,0);


\end{scope}

\end{scope}

    \end{tikzpicture}
  \caption{Dual triangulations to Hex and Y boards.}
  \label{fig:dualtriangle}
\end{figure}

  We present games of Hex and Y with boards being triangulations of a disk, with labels on the boundary vertices, such that players color vertices. To obtain this view from the classical one, one takes a dual triangulation to the Hex or the Y board, and then one labels the sides accordingly; see Figure~\ref{fig:dualtriangle}.

  \textbf{Board $B$:} Let $B$ be a triangulation of a disk. A \emph{path} between two vertices, $a$ and $b$ in $B$ is a sequence of vertices $a=v_0,v_1, \ldots, v_k=b$ such that any two consecutive vertices are connected by an edge. The boundary cycle of $B$ splits into a concatenation of $k$ paths $\alpha_1, \ldots, \alpha_k$. That is, the end vertex of path $\alpha_i$ is the beginning vertex of path $\alpha_{i+1}$ (indices taken mod $k$). We call paths $\alpha_i$ the \emph{sides} of $B$. To avoid certain degenerate cases, we assume that every $\alpha_i$ contains at least one edge (and thus at least two vertices).


  \textbf{Games of Hex and Y:} Players take turns to color the vertices of the triangulation of board $B$ with red and blue. One can color only non-colored vertices, and once a vertex is colored, it stays colored until the end of the game. The games differ in the shape of the board and the winning condition. A subset of vertices $A$ of $B$ is \emph{connected} if for any two vertices in $A$ there is a path in $A$ joining these two vertices. A \emph{chain} is a connected subset of vertices having the same color. 

  \textbf{Game of Hex:} 
  The board has $4$ boundary paths, called $R_1, B_1, R_2, B_2$. The goal for the red player is to create a red chain with at least one vertex in $R_1$ and at least one vertex in $R_2$. The goal for the blue player is to create a blue chain with at least one vertex in $B_1$ and at least one vertex in $B_2$. 

  \textbf{Game of Y:}
  The board has $3$ boundary paths $l_1, l_2,l_3$. The goal for either player is to form a chain of their color, containing at least one vertex in $l_1, l_2$, and $l_3$.

\section{Equivalence of Hex Theorem and Y Theorem}

\begin{figure}[!h]
  \centering
    \begin{tikzpicture}[scale=0.25]

      \definecolor{vlgray}{RGB}{230,230,230}
      \definecolor{red}{RGB}{255,000,000}

      \definecolor{blue}{RGB}{000,000,255}

\begin{scope}[shift={(-12,0)}]


\foreach \j in {0,...,5}{
  \foreach \i in {0,...,5}{
    
      \begin{scope}[shift={(2*\i-\j,\j*1.75)}]

      \draw[fill=black, black] (1,0.5)   circle [radius=0.1];
        \end{scope}
    }
}

\foreach \j in {0,...,5}{

      \begin{scope}[shift={(-\j,\j*1.75)}]
          \draw (1,0.5)--(11,0.5);

        \end{scope}

            \begin{scope}[shift={(2*\j,0)}]
          \draw (1,0.5)--(-4,9.25);
          
        \end{scope}
}

\foreach \i in {0,...,5}{ 

\begin{scope}[shift={(10-2*\i,0)}]

          \draw (1,0.5)--(1+ \i,0.5+ 1.75*\i);

\end{scope}

}

\foreach \i in {0,...,5}{ 

\begin{scope}[shift={(-5+2*\i,8.75)}]

          \draw (1,0.5)--(1-\i,0.5- 1.75*\i);

\end{scope}

}

\draw[ red, opacity=0.4, ultra thick] (1,0.5)--(11,0.5);

\draw[ red, opacity=0.4, ultra thick] (-4,9.25)--(6,9.25);

\draw[ blue, opacity=0.4, ultra thick] (1,0.5)--(-4,9.25);

\draw[ blue, opacity=0.4, ultra thick] (6,9.25)--(11,0.5);

\node [blue]  at (10.25,5.25)  {$B_2$};

\node [red]  at (6,-0.5)  {$R_1$};

\node [red]  at (1,10.5)  {$R_2$};

\node [blue]  at (-3.1,4.35)  {$B_1$};

\end{scope}

\begin{scope}[shift={(12,0)}]


\node [blue,left ]  at (-5.5,3)  {$b_0$};

\node [red,above ]  at (1.4,13.5)  {$r_0$};

\node [green]  at (8.5,7.25)  {$l_2$};

\node [violet]  at (-3,11.5)  {$l_1$};
\node [orange]  at (4,-0.5)  {$l_3$};

\foreach \j in {0,...,5}{
  \foreach \i in {0,...,5}{
    
      \begin{scope}[shift={(2*\i-\j,\j*1.75)}]

      \draw[fill=black, black] (1,0.5)   circle [radius=0.1];
        \end{scope}
    }
}

\foreach \j in {0,...,5}{

      \begin{scope}[shift={(-\j,\j*1.75)}]
          \draw (1,0.5)--(11,0.5);

        \end{scope}

            \begin{scope}[shift={(2*\j,0)}]
          \draw (1,0.5)--(-4,9.25);
          
        \end{scope}
}

\foreach \i in {0,...,5}{ 

\begin{scope}[shift={(10-2*\i,0)}]

          \draw (1,0.5)--(1+ \i,0.5+ 1.75*\i);

\end{scope}

}

\foreach \i in {0,...,5}{ 

\begin{scope}[shift={(-5+2*\i,8.75)}]

          \draw (1,0.5)--(1-\i,0.5- 1.75*\i);

\end{scope}

}

\draw[ green, opacity=0.4, ultra thick] (1,13.25)--(6,9.25)--(11,0.5);

\draw[ violet, opacity=0.4, ultra thick] (1,13.25)--(-4,9.25)--(-5.5,3);

\draw[ orange, opacity=0.4, ultra thick] (-5.5,3)--(1,0.5)--(11,0.5);

 \foreach \i in {0,...,5}{ 
 \draw (-5.5,3)--(1-1*\i,0.5+1.75*\i);

}

 \draw[fill=blue, blue] (-5.5,3)   circle [radius=0.2];

  \foreach \i in {0,...,5}{ 
 \draw (1,13.25)--(-4+2*\i,9.25);

 \draw[fill=red, red] (1,13.25)   circle [radius=0.2];

}

\end{scope}

    \end{tikzpicture}
  \caption{Extending a Hex board to a position on a Y board.}
  \label{fig:hexbaordtoyboard}
\end{figure}
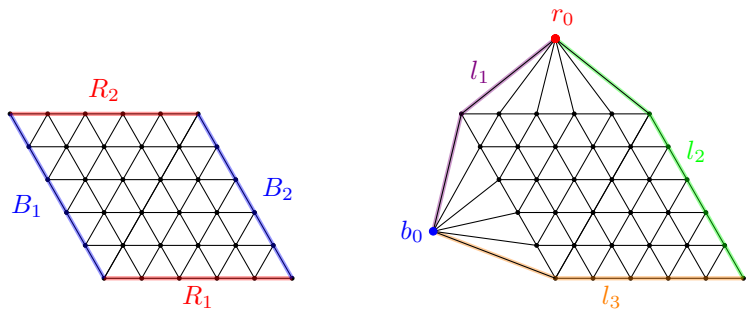

  \textbf{Y $\Rightarrow$ Hex.} This observation appears throughout the literature and is attributed to Schensted \cite{SchenTit}. Consider a Hex board $B$. Add one vertex  $r_0$ and join it by edges to all vertices of $R_2$. Add one vertex $b_0$ and join it by edges to all vertices of $B_1$. Color $r_0$ red, and color $b_0$ blue. Now we view the obtained board $B'$ as a $Y$ game, with $l_1=\{b_0\} \cup \{(B_1 \cap R_2)\}\cup \{ r_0\}$,  $l_2=B_2 \cup \{r_0\}$, and $l_3= R_1 \cup  \{b_0\}$. (Since  every side of  $B$ contains at least one edge, it is clear that $B'$ is a triangulation of a disk.) Boards $B$ and $B'$ are presented in Figure~\ref{fig:hexbaordtoyboard}.

  The reader easily sees that playing Y from this position onwards is the same as playing Hex on the original board. By Y Theorem, this particular position has a unique winner, thus so does the corresponding Hex game.

\begin{figure}[!h]
  \centering
    \begin{tikzpicture}[scale=0.45]

      \definecolor{vlgray}{RGB}{230,230,230}
      \definecolor{red}{RGB}{255,000,000}

      \definecolor{blue}{RGB}{000,000,255}

\begin{scope}[shift={(0,0)}]


\foreach \j in {0,...,5}{
  \foreach \i in {0,...,5}{
    
      \begin{scope}[shift={(2*\i-\j,\j*1.75)}]
        
      \draw[fill=black, black] (1,0.5)   circle [radius=0.1];
        \end{scope}
    }
}

\foreach \j in {0,...,5}{

      \begin{scope}[shift={(-\j,\j*1.75)}]
          \draw (1,0.5)--(11,0.5);

        \end{scope}

            \begin{scope}[shift={(2*\j,0)}]
          \draw (1,0.5)--(-4,9.25);
          
        \end{scope}
}

\foreach \i in {0,...,5}{ 

\begin{scope}[shift={(10-2*\i,0)}]

          \draw (1,0.5)--(1+ \i,0.5+ 1.75*\i);

\end{scope}

}

\foreach \i in {0,...,5}{ 

\begin{scope}[shift={(-5+2*\i,8.75)}]

          \draw (1,0.5)--(1-\i,0.5- 1.75*\i);

\end{scope}

}

\node [green]  at (10.25,5.25)  {$l_2=B_2$};

\node [orange]  at (6,-0.5)  {$l_3=R_1$};

\node [green]  at (1,10.5)  {$l_2'=R_2$};

\node [orange]  at (-3.1,4.35)  {$l_3'=B_1$};

\node [violet]  at (2.05,3.35)  {$l_1$};

\node []  at (7,1.75)  {$B$};
\node []  at (-1,6.35)  {$B'$};

\node [red]  at (2.1,6.9)  {$C'$};
\node [red]  at (6,4.6)  {$C$};

\draw[ orange, opacity=0.4, ultra thick] (1,0.5)--(11,0.5);

\draw[ green, opacity=0.4, ultra thick] (-4,9.25)--(6,9.25);

\draw[ violet, opacity=0.4, ultra thick] (1,0.5)--(6,9.25);

\draw[ orange, opacity=0.4, ultra thick] (1,0.5)--(-4,9.25);

\draw[ green, opacity=0.4, ultra thick] (6,9.25)--(11,0.5);

\draw [line width=4, red!40] (3,0.5)--(4,2.25)--(3,4)--(4,5.75)--(3,7.5)--(1,7.5)--(0,9.25);

\draw [line width=3, dashed, red] (4,5.75)--(6,5.75)--(7,4)--(9,4);

 \draw [line width=3, dashed, red]  (3,0.5)--(4,2.25)--(3,4)--(4,5.75);

 \draw[fill=red, red] (3,0.5)  circle [radius=0.2];

 \draw[fill=red, red] (4,2.25)  circle [radius=0.2];
 \draw[fill=red, red] (3,4)  circle [radius=0.2];
 \draw[fill=red, red] (4,5.75)  circle [radius=0.2];
 \draw[fill=red, red] (3,7.5) circle [radius=0.2];
 \draw[fill=red, red] (1,7.5)  circle [radius=0.2];
\draw[fill=red, red] (0,9.25)  circle [radius=0.2];

 \draw[fill=red, red] (6,5.75) circle [radius=0.2];
  \draw[fill=red, red] (7,4) circle [radius=0.2];
  \draw[fill=red, red] (9,4) circle [radius=0.2];

\draw [thick,->] (5,11) [out=10, in=120] to (8,9.5);
\node []  at (7,11.5)  {$p$};

\end{scope}

    \end{tikzpicture}
  \caption{Doubling the Y board, and folding it back onto itself. Chain $C'$ is pink. Chain $C$ is dashed red.}
  \label{fig:folding}
\end{figure}
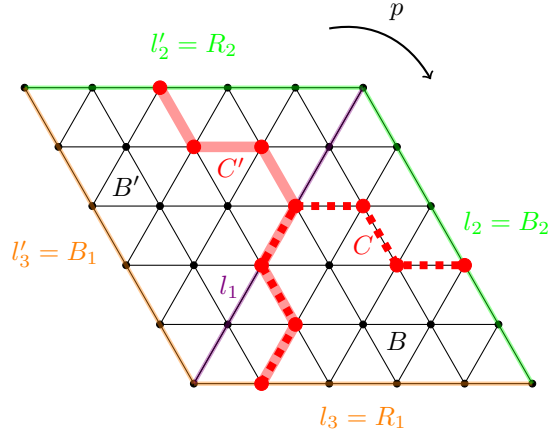

  \textbf{Hex $\Rightarrow$ Y.} First we prove that there exists at least one winner. Consider a Y board $B$ with sides $l_1, l_2$, and $l_3$. Assume that every vertex of $B$ is colored. Take another copy $B'$ of $B$ (together with the coloring) and glue it to $B$ along the identity map on $l_1$. Since $l_1$ has at least one edge, the resulting space $B \cup_{l_1} B'$ is a triangulation of a disk with $4$ boundary paths. Let $l_2'$ and $l_3'$ denote the sides of $B \cup_{l_1} B'$ opposite to $l_3$ and $l_2$ respectively. Treat $B \cup_{l_1} B'$ as a Hex board, where $l_3=R_1, l_2=B_2, l_3'=B_1$ and $l_2'=R_2$. The board $B \cup_{l_1} B'$ is shown in Figure~\ref{fig:folding}. By Hex Theorem there is either a red chain joining $R_1$ and $R_2$ or a blue chain joining $B_1$ and $B_2$. Assume that it is a red chain and call it $C'$. The case of the blue chain is done analogously.

  Note that since $l_1$ is a disconnecting subset of $ B \cup_{l_1} B'$, chain $C'$ has to have at least one vertex in $l_1$.  Now let $C$ be a subset of $B$ which is the image of $C'$ under the map $p \colon B \cup_{l_1} B' \to B$ which folds $B'$ onto $B$. Note that $p$ is simplicial (i.e.,\ if two vertices are connected by an edge, then so are their images), color-preserving and it restricts to the identity map on $B$. Therefore $C=p(C')$ is a chain as well. Since $C'$ has a vertex in $R_2= l_2'$, and since $p(l_2')=l_2$, we get that $C$ has vertices on all three sides of $B$, thus ensuring a win for red. 

\begin{figure}[!h]
  \centering
    \begin{tikzpicture}[scale=0.45]

      \definecolor{vlgray}{RGB}{230,230,230}
      \definecolor{red}{RGB}{255,000,000}

      \definecolor{blue}{RGB}{000,000,255}

\begin{scope}[shift={(0,0)}]


\foreach \j in {0,...,5}{
  \foreach \i in {0,...,5}{
    
      \begin{scope}[shift={(2*\i-\j,\j*1.75)}]
        
      \draw[fill=black, black] (1,0.5)   circle [radius=0.1];
        \end{scope}
    }
}

\foreach \j in {0,...,5}{

      \begin{scope}[shift={(-\j,\j*1.75)}]
          \draw (1,0.5)--(11,0.5);

        \end{scope}

            \begin{scope}[shift={(2*\j,0)}]
          \draw (1,0.5)--(-4,9.25);
          
        \end{scope}
}

\foreach \i in {0,...,5}{ 

\begin{scope}[shift={(10-2*\i,0)}]

          \draw (1,0.5)--(1+ \i,0.5+ 1.75*\i);

\end{scope}

}

\foreach \i in {0,...,5}{ 

\begin{scope}[shift={(-5+2*\i,8.75)}]

          \draw (1,0.5)--(1-\i,0.5- 1.75*\i);

\end{scope}

}

\node [green]  at (10.25,5.25)  {$l_2=B_2$};

\node [orange]  at (6,-0.5)  {$l_3=R_1$};

\node [green]  at (1,10.5)  {$l_2'=R_2$};

\node [orange]  at (-3.1,4.35)  {$l_3'=B_1$};

\node [violet]  at (2.05,3.35)  {$l_1$};

\node []  at (7,1.75)  {$B$};
\node []  at (-1,6.35)  {$B'$};

\node [red]  at (2.05,6.9)  {$s(C)$};
\node [red]  at (6,4.6)  {$C$};

\draw[ orange, opacity=0.4, ultra thick] (1,0.5)--(11,0.5);

\draw[ green, opacity=0.4, ultra thick] (-4,9.25)--(6,9.25);

\draw[ violet, opacity=0.4, ultra thick] (1,0.5)--(6,9.25);

\draw[ orange, opacity=0.4, ultra thick] (1,0.5)--(-4,9.25);

\draw[ green, opacity=0.4, ultra thick] (6,9.25)--(11,0.5);

\draw [line width=4, red!40]  (5,0.5)--(6,2.25)--(7,4)--(9,4);
\draw [line width=4, red!40]  (7,4)--(6,5.75)--(4,5.75);

 \draw[fill=red, red] (5,0.5)  circle [radius=0.2];

 \draw[fill=red, red] (6,2.25)  circle [radius=0.2];
 \draw[fill=red, red] (7,4)  circle [radius=0.2];
 \draw[fill=red, red] (9,4)  circle [radius=0.2];
 \draw[fill=red, red] (6,5.75) circle [radius=0.2];
 \draw[fill=red, red] (4,5.75)  circle [radius=0.2];

\draw [line width=3, dashed, red] (4,5.75)--(3,7.5);

\begin{scope}[shift={(-6,3.5)}]

\draw [line width=3, dashed, red] (5,0.5)--(6,2.25)--(7,4)--(9,4);
\draw [line width=3, dashed, red] (7,4)--(6,5.75);

 \draw[fill=red, red] (5,0.5)  circle [radius=0.2];

 \draw[fill=red, red] (6,2.25)  circle [radius=0.2];
 \draw[fill=red, red] (7,4)  circle [radius=0.2];
 \draw[fill=red, red] (9,4)  circle [radius=0.2];
 \draw[fill=red, red] (6,5.75) circle [radius=0.2];

\end{scope}

\draw [thick,<->] (5,11) [out=10, in=120] to (8,9.5);
\node []  at (7,11.5)  {$s$};

\end{scope}

    \end{tikzpicture}
  \caption{A winner on a Y board leads to a winner on a Hex board. Chain $C$ is pink. Chain $s(C)$ is dashed red.}
  \label{fig:winnertowinner}
\end{figure}
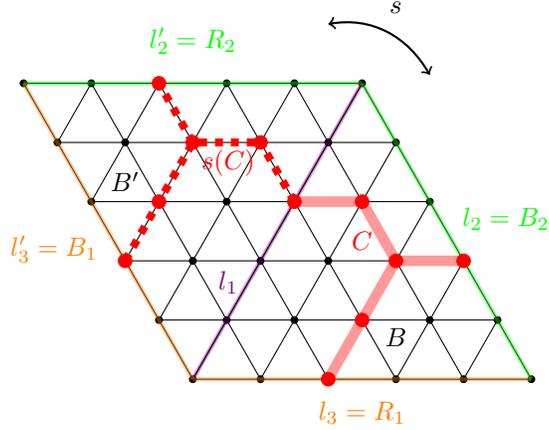

  Now we show that there is a unique winner. Let $B , B'$ and $B \cup_{l_1} B'$ be as above and let $s \colon B \cup_{l_1} B' \to B \cup_{l_1} B'$ be the reflection across $l_1$. Now suppose we have a chain $C$ connecting all three sides of $B$. Then the union $C \cup s(C)$ is a chain connecting all four sides of $B \cup_{l_1} B'$; see Figure~\ref{fig:winnertowinner}. Thus, having two winners of the $Y$ game on $B$ leads to two winners of the corresponding Hex game on $B \cup_{l_1} B'  $.

\section{Proof of Hex Theorem}

  We present essentially the same proof as Gale \cite{Gale}, with one nuance: rather than defining a fixed-point free map of a disk to itself, we use our generalized representation of a board to directly construct a retraction from a disk onto its boundary circle.

  Consider a Hex board $B$, with sides $R_1, B_1, R_2, B_2$. Add four vertices $r^{-}$, $r^{+}$, $b^{-}$, $b^{+}$ and connect $r^{-}$ to both $b^{-}$ and $b^{+}$, and $r^{+}$ to both $b^{-}$ and $b^{+}$. Then connect every vertex of $R_1$ to $r^{-}$, every vertex of $R_2$ to $r^{+}$, every vertex of $B_1$ to $b^{-}$, and every vertex of $B_2$ to $b^{+}$. Thus we get a triangulation of a slightly larger disk $D$, which is presented in Figure~\ref{fig:proofofhex}. The original board and the game are expressed in terms of the new board as follows. The boundary vertices $r^{-}$ and $r^{+}$ are (from the beginning of the game) colored red, and the vertices $b^{-}$ and $b^{+}$ are colored blue. Hex Theorem for $B$ is equivalent to the following statement for $D$: there is either a red chain form $r^{-}$ to $r^{+}$, or a blue chain from $b^{-}$ to $b^{+}$, but not both.

  \begin{proof} First we prove that there is at least one winner. Assume that every vertex of $D$ is colored. Let $V^{+}$ denote the subset of all red vertices of $D$ which are connected by a red chain to $r^{+}$. Let $V^{-}$ be defined as all the red vertices that are not in $r^{-}$. Define subsets of $W^{+}$ and $W^{-}$ for the blue vertices analogously. By definition, there is no edge between a vertex in $V^{+}$ and a vertex of $V^{-}$, and the same holds for $W^{+}$ and $W^{-}$. Note that also by definition we have that $r^{+}$ is in $V^{+}$ and $b^{+}$ is in $W^{+}$.

  Suppose by contrary that there is no winner. This implies that $r^{-}$ is not in $V^{+}$ and thus it is in $V^{-}$. For the same reason we have that $b^{-}$ is in $W^{-}$. Let $S = (r^{-}, b^{-}, r^{+}, b^{+})$ denote the boundary cycle of $D$. We define a map $ D \to S$  on vertices by sending the above subsets as follows:
  \begin{align*}
                V^{-}  &\mapsto r^{-}, \\
                V^{+}  &\mapsto r^{+}, \\
                W^{-}  &\mapsto b^{-}, \\
                W^{+}  &\mapsto b^{+}. 
  \end{align*}

  Observe that this assignment gives a simplicial map. Therefore it induces a continuous map $ D \to S$ with the topology coming from the standard Euclidean metric on every triangle. One easily sees that this map is the identity on $S$, which gives a retraction of the disk onto its boundary circle, a contradiction \cite[Corollary~2.15]{Hat}.

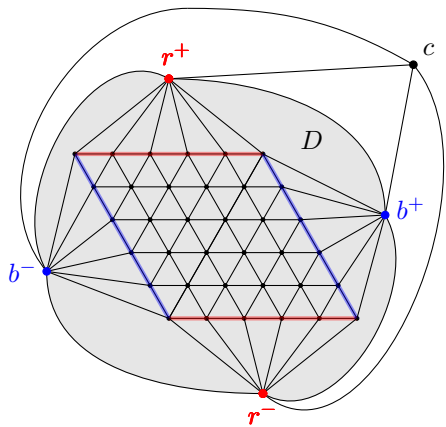
\begin{figure}[!h]
  \centering
    \begin{tikzpicture}[scale=0.25]

      \definecolor{vlgray}{RGB}{230,230,230}
      \definecolor{red}{RGB}{255,000,000}

      \definecolor{blue}{RGB}{000,000,255}

\begin{scope}[shift={(0,0)}]

\draw [fill=vlgray, vlgray] (12.5,6) [out=300,in=330] to (6,-3.5) [in=270,out=180] to (-5.5,3) [in=150,out=120] to (1,13.25) [in=90,out=0] to (12.5,6);


\foreach \j in {0,...,5}{
  \foreach \i in {0,...,5}{
    
      \begin{scope}[shift={(2*\i-\j,\j*1.75)}]

      \draw[fill=black, black] (1,0.5)   circle [radius=0.1];
        \end{scope}
    }
}

\foreach \j in {0,...,5}{

      \begin{scope}[shift={(-\j,\j*1.75)}]
          \draw (1,0.5)--(11,0.5);

        \end{scope}

            \begin{scope}[shift={(2*\j,0)}]
          \draw (1,0.5)--(-4,9.25);
          
        \end{scope}
}

\foreach \i in {0,...,5}{ 

\begin{scope}[shift={(10-2*\i,0)}]

          \draw (1,0.5)--(1+ \i,0.5+ 1.75*\i);

\end{scope}

}

\foreach \i in {0,...,5}{ 

\begin{scope}[shift={(-5+2*\i,8.75)}]

          \draw (1,0.5)--(1-\i,0.5- 1.75*\i);

\end{scope}

}

\draw[ red, opacity=0.4, ultra thick] (1,0.5)--(11,0.5);

\draw[ red, opacity=0.4, ultra thick] (-4,9.25)--(6,9.25);

\draw[ blue, opacity=0.4, ultra thick] (1,0.5)--(-4,9.25);

\draw[ blue, opacity=0.4, ultra thick] (6,9.25)--(11,0.5);

\draw (-5.5,3) [in=180,out=130] to   (2,17);

\draw (2,17) [in=150,out=0] to  (14,14);

\draw (1,13.25)  to (14,14);

\draw (6,-3.5) [in=310,out=320] to (14,14);
\draw  (12.5,6)  to (14,14);

\draw (-5.5,3) [in=150,out=120] to (1,13.25);
\draw (-5.5,3) [in=180,out=270] to (6,-3.5);
\draw (1,13.25) [in=90,out=0] to (12.5,6);

\draw (6,-3.5) [in=300,out=330] to (12.5,6);

 \foreach \i in {0,...,5}{ 
 \draw (-5.5,3)--(1-1*\i,0.5+1.75*\i);

}

 \draw[fill=blue, blue] (-5.5,3)   circle [radius=0.2];

\node [blue,left ]  at (-5.5,3)  {$b^{-}$};

\begin{scope}[shift={(10,0)}]

 \foreach \i in {0,...,5}{ 
 \draw (2.5,6)--(1-1*\i,0.5+1.75*\i);

}

 \draw[fill=blue, blue] (2.5,6)   circle [radius=0.2];
\node [blue,right ]  at (2.6,6.35)  {$b^{+}$};

\end{scope}

  \foreach \i in {0,...,5}{ 
 \draw (1,13.25)--(-4+2*\i,9.25);

 \draw[fill=red, red] (1,13.25)   circle [radius=0.2];
\node [red,above ]  at (1.4,13.5)  {$r^{+}$};

}

\begin{scope}[shift={(5,9.75)},yscale=-1]

  \foreach \i in {0,...,5}{ 
  \draw (1,13.25)--(-4+2*\i,9.25);

  \draw[fill=red, red] (1,13.25)   circle [radius=0.2];

\node [red,below ]  at (1,13.25)  {$r^{-}$};

}
\end{scope}

 \draw[fill=black] (14,14)   circle [radius=0.2];
\node [above right]  at (14,14)  {$c$};

\node [above right]  at (7.5,9)  {$D$};

\end{scope}

    \end{tikzpicture}
  \caption{A $2$--disk $D$ enclosing a Hex board.}
  \label{fig:proofofhex}
\end{figure}

  It remains to show that there cannot be two winners. Assume the contrary. Thus we have chains of respective color joining  $r^{-}$ to $r^{+}$ and $b^{-}$ to $b^{+}$, which by definition lie inside of $D$ and do not intersect each other. By taking such chains with minimal number of vertices we can assume that they are embedded paths. Embed $D$ into the plane, add one more vertex $c$ on the exterior of $D$, and connect all vertices $r^{-}, r^{+}, b^{-}, b^{+}$ to $c$ by edges such that they do not intersect one another (the reader easily sees that it can always be done, we present one such configuration in Figure~\ref{fig:proofofhex}). This gives an embedding of the topological (i.e.,\ isomorphic after forgetting vertices of degree $2$) complete graph $K_5$ into the plane which is impossible \cite{Kuratowski}.
\end{proof}

\end{document}